\documentclass{elsart}

\usepackage{epsfig}

\usepackage{amsmath}
\usepackage{amsfonts}
\usepackage{amscd}
\usepackage{amssymb}
\usepackage{diagrams}

\def\CC{{\mathbb C}}
\def\ZZ{{\mathbb Z}}
\def\NN{{\mathbb N}}
\def\RR{{\mathbb R}}
\def\PP{{\mathbb P}}
\def\TT{{\mathcal T}}
\def\O{{\mathcal O}}
\def\F{{\mathcal F}}
\def\res{{\rm res}}
\def\Hom{{\rm Hom}}
\def\floor#1{\lfloor #1 \rfloor}
\def\ceil#1{\lceil #1 \rceil}

\begin{document}

\begin{frontmatter}



\title{Exact Matrix Formula for the Unmixed Resultant in Three Variables}
\author{Amit Khetan}
\ead{khetan@math.umass.edu}
\address{Department of Mathematics, University of Massachusetts, Amherst, MA, USA}




\abstract

We give the first exact determinantal formula for the
resultant of an unmixed sparse system of four Laurent polynomials in
three variables with arbitrary support.  This follows earlier work by
the author on exact formulas for bivariate systems and also uses the
exterior algebra techniques of Eisenbud and Schreyer. Along the way we
will prove an interesting new vanishing theorem for the sheaf
cohomology of divisors on toric varieties. This will also allow us to
describe some supports in four or more variables for which
determinantal formulas for the resultant exist.

\endabstract

\end{frontmatter}

\section{Introduction}
The resultant of $n+1$ polynomials $f_1, \dots, f_{n+1}$ in $n$
variables is a single polynomial in the coefficients of the $f_i$
which vanishes when the $f_i$ have a common root. The resultant
can therefore be used to eliminate $n$ variables from $n+1$
equations. Originally resultants were defined for generic
polynomials of fixed total degrees. More recently a {\em sparse
resultant} has been defined which exploits the monomial structure
of the given polynomials. The foundational work was laid by
Kapranov, Sturmfels, and Zelevinsky \cite{KSZ}. Sparse resultants
are discussed in depth in the book \cite{GKZ}.

Formally, let $f_1, f_2, \dots, f_{n+1} \in \CC[x_1, x_1^{-1}, \dots,
x_n, x_n^{-1}]$ be polynomials with the same Newton polytope $Q$. Let
$A = Q \cap \ZZ^n$. We will assume that $A$ affinely generates
$\ZZ^n$.

We can write:
$$f_i = \sum_{\alpha \in A} C_{i\alpha} x^\alpha $$
We will treat the coefficients $C_{i\alpha}$ as
independent variables throughout.

\begin{defn} The $A${\em -resultant} $\res_A(f_1, \dots, f_{n+1})$ is the irreducible polynomial in the ring $\ZZ[C_{i\alpha}]$, unique up to sign,
which vanishes whenever $f_1, \dots, f_n$ have a common root in
$(\CC^\ast)^n$.
\end{defn}

The problem of finding explicit formulas for resultants, and their
cousins the discriminants, dates back to the nineteenth century with
the work of Cayley, Sylvester, B\'ezout and others. With the recent
increase in computing power there has been a renewed interest in
computing resultants and new applications in fields such as computer
graphics, machine vision, robotic inverse kinematics, and molecular
structure \cite{MC1}, \cite{MC2}, \cite{EPhD}.

Even in very small examples, the resultant can have millions of terms.
Therefore most authors have looked for a more compact representation.
A determinantal formula, following the classical formulas of Sylvester
and B\'ezout, writes the resultant as the determinant of a matrix
whose entries are easily computable polynomials of low degree. In the
dense case, when all the polynomials have the same degree,
determinantal formulas are known when $n=1, 2$, or $3$ and for a very
few cases in more variables. In the sparse case, $n=1$ is the same as
the dense case and there are the classical Sylvester and B\'ezout
formulas, determinantal formulas for $n=2$ were found by the author in
\cite{khe}. This paper gives a new exact formula when $n=3$.

Given any lattice polytope $Q$, let $D_1, \dots, D_s$ denote the
facets (codimension 1 faces) of $Q$. Given a subset $I = \{i_1, \dots,
i_k\}$ of $\{1, \dots, s \}$, let $D_I =\{ D_{i_1}, \dots, D_{i_k} \}$
be the corresponding subset of facets. Let $\overline{D_I}$ be the set
of facets of $Q$ not in $D_I$.  $Q - D_I$ will refer to the set of all
points in $Q$ but not on any facet on $D_I$. More generally, $kQ -
D_I$ is the set of all points in the Minkowski sum of $k$ copies of
$Q$ but not on any of the facets corresponding to $D_I$. Finally,
given a set $S \subset \RR^n$ let $l(S) = S \cap \ZZ^n$ be the set of
{\em lattice} points in $S$.  The main theorem is as follows:

\begin{thm} \label{thm1}
Let $f_1, f_2, f_3, f_4 \in \CC[x_1, x_2, x_3, x_1^{-1}, x_2^{-1},
x_3^{-1}]$ be four polynomials with common Newton polytope $Q \subset
\RR^3$. Suppose $A = Q \cap \ZZ^3$ affinely generates $\ZZ^3$. Pick a
proper collection of the facets of $Q$, $D_{I} = (D_{i_1}, \dots,
D_{i_k})$, such that the union of the facets in $D_I$ is homeomorphic
to a disk. There is a determinantal formula for the resultant
$\res_A(f_1, f_2, f_3, f_4)$ of the following block form:
$$
\begin{pmatrix}
 B &  L \\ \tilde L & 0.
 \end{pmatrix}
$$
\noindent
The rows of $B$ and $L$ are indexed by the points in $l(2Q-D_I)$. The
columns of $B$ and $\tilde{L}$ are indexed by $l(2Q -
\overline{D_I})$.  The rows of $\tilde{L}$ are indexed by four copies
of $l(Q - \overline{D_I})$, and the columns of $L$ are indexed by four
copies of $l(Q - D_I)$.

The entries of $B$ are of B{\'e}zout type and are polynomials of
degree $4$ in the coefficients $C_{i\alpha}$.  The entries of $L$ and
$\tilde{L}$ are of Sylvester type, thus linear in the $C_{i\alpha}$.

\end{thm}

We will see how the entries of $B$ can be filled in using a free
resolution over an exterior algebra. Both the proof and the
construction are based on techniques developed by Eisenbud and
Schreyer, which have been adapted for sparse resultants (toric
varieties).

The paper is organized as follows. Section \ref{s:setup} discusses the
background on toric varietes, exterior algebras, and the Tate
resolution of Eisenbud-Schreyer. Section \ref{s:proof} uses these
techniques along with some sheaf cohomology vanishing results to prove
Theorem \ref{thm1}. In particular Section \ref{s:proof} contains a new
vanishing result for certain divisors on any projective toric
variety. Section \ref{s:construct} shows how to actually construct the
resultant matrix and gives some examples. Finally, Section
\ref{s:ehrhart} gives a different combinatorial perspective on the
resultant matrix in terms of the Ehrhart polynomial and analyzes the
size of the resultant matrix.

\section{Notation and Background}

\label{s:setup}

\subsection{Toric varieties and Chow forms}

\label{ss:toric}

Given a polytope $Q \subset \RR^n$ and associated $A = Q \cap \ZZ^n$,
let $N = |A|$. The {\em toric variety} $X_A \subset \PP^{N-1}$ is
defined as the algebraic closure of the set $(x^{\alpha_1}: \dots:
x^{\alpha_N})$ where $\alpha_i$ ranges over the elements of $A$ and $x
\in (\CC^\ast)^n$. It has dimension $n$. In terms of $X_A$, the
polynomials $f_i$ are hyperplane sections. The system $(f_1, f_2,
\dots, f_{n+1})$ defines a codimension $n+1$ plane. The set of all
codimension $n+1$ planes meeting $X_A$ defines a hypersurface in the
Grasmannian $G(n+1, N)$. The $A$-resultant is identified with the
equation of this hypersurface, also called the {\em Chow form} of
$X_A$.

\begin{prop}
The resultant $\res_A(f_1, \dots, f_{n+1}) = 0$ if and only if the
$f_i$ have a common intersection on $X_A$.
\end{prop}

Let $\Sigma_Q$ be the normal fan of $Q$ with $\Sigma_Q(1) = \{\eta_1,
\dots, \eta_s\}$ the inner normals to the facets. There is an
associated normal toric variety $X_{\Sigma_Q}$ (see \cite[Chapter
1]{Ful}). Assuming $A$ affinely spans $\ZZ^n$, $X_{\Sigma_Q}$ is the
normalization of $X_A$. This is essentially Proposition 4.9 in Chapter
5 of \cite{GKZ}.  The results below are standard and can be found in
\cite{Ful}.

\begin{prop}
\label{p:toricdivisors} The $\eta_i$ are in 1-1 correspondence
with the torus invariant prime Weil divisors on $X_{\Sigma_Q}$.
Let $D_i$ denote the divisor corresponding to $\eta_i$, and
$\O(D_i)$ the corresponding rank 1 reflexive sheaf on
$X_{\Sigma_Q}$.
\end{prop}

In the introduction, and statement of Theorem \ref{thm1}, $D_i$
denoted a facet of $Q$. This facet will be identified with the
corresponding prime divisor, also denoted $D_i$, as defined above.

Given a general divisor $D = \sum a_i D_i$ on $X_{\Sigma_Q}$, we will
denote by $\O_{X_A}(D)$ or, when there is no confusion just $\O(D)$,
the push-forward of the sheaf $\O_{X_{\Sigma_Q}}(D)$ onto $X_A$ via
the normalization map. The linear equivalence classes of divisors are
computed by the following exact sequence:
$$ \begin{CD} 0 @>>> \ZZ^n @>div>> \ZZ^s @>{[\cdot]}>> Cl_X @>>> 0, \end{CD}$$
\noindent where $div(u) = (\langle u, \eta_1 \rangle, \dots, \langle u, 
\eta_s, \rangle)$ and $Cl_X$ is the cokernel of this map. Given a divisor $D
\in \ZZ^s$ we let $[D]$ be the image of $D$ in $Cl_X$.

There is a nice combinatorial description of the global sections
$H^0(X_A, \O(D))$.  A divisor $D = \sum a_i D_i$ determines a convex
polytope $P_D = \{ m \in \RR^n \ : \ \langle m, \eta_i \rangle \geq
-a_i \}$. For any polytope $P$, let $S_P$ denote the $\CC$ vector
space with basis the lattice points in $P$, i.e., $S_P = \CC\{P \cap
\ZZ^n\}$.

\begin{prop} \label{p:h0prop} $$H^0(X_A, \O(D)) \cong S_{P_D}.$$
\end{prop}

If we start with a polytope $Q$, then it determines an ample divisor
on the toric variety $X_{\Sigma_Q}$. Write:

\[ Q = \{ m \in \RR^n \: \ \langle m, \eta_i \rangle \geq -a_i, \ i = 1, \dots, s \}, \]

\noindent for some $a_1, \dots, a_s \in \ZZ.$ Let $D_Q = \sum a_i D_i$
be the corresponding divisor. If $X_A$ is the (possibly non-normal)
toric variety above defined by the lattice points in $Q$, then the
push-forward of $D_Q$ yields the very ample divisor corresponding to
the embedding of $X_A$ into $\PP^{N-1}$. On $X_{\Sigma_Q}$, $D_Q$ will
always be ample but not necessarily very ample. One final useful fact
is that the  sheaf $\mathcal{O}(-\sum_{i=1}^s D_i)$ is the
canonical sheaf on the Cohen-Macaulay variety $X_{\Sigma_Q}$. This
will be needed when we apply Serre duality below.

\subsection{Exterior algebra and the Tate resolution}

\label{ss:tate}

Eisenbud and Schreyer \cite{ES} have developed some powerful new
machinery to compute Chow forms using resolutions over an exterior
algebra. Suppose $X \subset \PP^{N-1}$ is a variety of dimension
$n$. We are interested in the finding the Chow form of $X$.

The ambient projective space $\mathbf{P} = \PP^{N-1}$ has the graded
coordinate ring $R = \CC[X_1, \dots, X_N]$. If we let $W$ be the $\CC$
vector space spanned by the $X_i$, (identified with the degree 1 part
of $R$), then $\mathbf{P}$ is the projectivization $\PP(W)$. The ring
$R$ can also be identified with the symmetric algebra ${\rm Sym}(W)$.

Now let $V = W^\ast$, the dual vector space, with a corresponding dual
basis $e_1, \dots e_N$.  We will consider the {\em exterior algebra}
$E = \bigwedge V$, also a graded algebra where the generators $e_i$
have degree $-1$.  We will use the standard notation $E(k)$ to refer
the rank 1 free $E$-module generated in degree $-k$.

For any coherent sheaf $\F$ on $\mathbf{P}$, there is an
associated exact complex of graded free $E$-modules, called the {\em
Tate resolution}, denoted $T(\F)$. The terms of
$T(\F)$ can be written in terms of the vector spaces of sheaf
cohomology of twists of $\F$. Namely, we
have:

\begin{equation} \label{tate}
T^{e}(\F) = \oplus [H^j(\F(e-j)) \otimes_{\CC} E(j-e)]
\end{equation}

Here $e$ is any positive integer. In particular, this complex is
infinite in both directions, although the terms themselves are finite
dimensional free $E$-modules.

Now suppose that $\mathcal{F}$ is supported on $X$. Recall that the
Chow form of $X$, also called the $X$-resultant and denoted $\res_X$,
is the defining equation of the set of codimension $n+1$-planes
meeting $X$.  Such a plane is specified by a $n+1$ dimensional
subspace $W_f = \CC \{ f_1, \dots, f_{n+1} \} \subset W$. Let
$\mathbf{G}$ be the Grasmannian of codimension $n+1$-planes on
$\mathbf{P}$. Let $ \mathcal{T}$ be the tautological bundle on
$\mathbf{G}$, that is to say the fiber at the point corresponding to
$f$ is just $W_f$. There is a functor, $U_{n+1}$ from free $E$-modules
to vector bundles on $\mathbf{G}$ which sends $E(p)$ to $\wedge^p
\mathcal{T}$.

This functor when applied to the Tate resolution gives a {\em finite}
complex of vector bundles on $\mathbf{G}$, $U_{n+1}(T(\F))$ that is
fiberwise a finite complex of $\CC$ vector spaces.

\begin{thm}
$$\det(U_{n+1}(T(\F))) = {\res_X}^{{\rm rank}(\F)}$$
\end{thm}

This is a determinant of a complex, which in general can be computed
as a certain alternating product of determinants. We will be most
interested in the special case where the complex in question has only
two terms:
$$ \begin{CD} 0 @>>> A @>\Psi>>B @>>> 0 \end{CD} $$
In this case, the determinant of the complex is just the determinant
of the matrix of the map $\Psi$. Sheaves whose Tate resolutions yield
such two term complexes for the Chow form are called {\em weakly
Ulrich}.  Determinantal formulas for the resultant correspond to
finding a weakly Ulrich sheaf of rank $1$ on the toric variety $X_A$.

Let $M = \oplus_{i \in \NN} H^0(\F(i))$. This is a graded
$R$-module. The {\em linear strand} of the Tate resolution is the
subcomplex defined by the terms $M_e \otimes E(-e)$. The maps in the linear
strand are completely canonical:
\begin{align*}
\phi_e \ :  \ M_{e} \otimes E(-e) &\to M_{e+1} \otimes E(-e-1) \\
m \otimes 1 &\mapsto \sum_{i=1}^N m \cdot X_i \otimes e_i
\end{align*}
An extremely important fact is that for large enough $e$, anything
larger than the {\em regularity} of $M$, all the higher
cohomology vanishes and only the linear strand remains.
For a definition and discussion on regularity see \cite{BB}.

This suggests an algorithm to compute terms of the Tate resolution:

\begin{enumerate} \label{tatealg}
\item Given $\F$ compute $M$.
\item Pick $e = {\rm reg}(M)+1$ and compute $\phi_e$.
\item Start computing a free resolution of $\phi_e$ over $E$.
\end{enumerate}

\noindent Note: As a consequence we can read off the cohomology of
twists of $\F$ as graded pieces of this resolution. As Eisenbud,
Schreyer, and Fl{\o}ystad \cite{EFS} point out, in many cases this
is the most efficient known way to compute sheaf cohomology.

\section{Proof of Theorem \ref{thm1}}

\label{s:proof}

Suppose we are given $f_1, f_2, f_3, f_4$ with common Newton polytope
$Q \subset \RR^3$.  To apply the exterior algebra construction we
take $W = S_Q$, the $\CC$ vector space with basis the lattice
points in $Q$, and $V = S_Q^{\ast}$.  The corresponding projective
space is $\mathbf{P} = \PP(W) \cong \mathbf{P}^{N-1}$, and the
exterior algebra is $E = \bigwedge V$. Let $y_1, \dots, y_N$
denote the basis of $S_Q$ and $e_1, \dots, e_N$ the corresponding
dual basis of $E$.

We now show how Theorem \ref{thm1} reduces to showing that an
appropriate push-forward of a Weil divisor class onto $X_A$ is a
weakly Ulrich sheaf. This will require proving that certain
cohomology groups vanish.

Let $I \subset \{1, \dots, s\}$, thought of as a subset of the
facets. Let $D_I = \sum_{i \in I} D_i$ and $\overline D_I = \sum_{i
\notin I} D_i$ be formal sums of the corresponding divisors. The
sheaves we will be interested in are of the form $kD_Q - D_I$ where $k
\in \ZZ$.

As in the statement of Theorem \ref{thm1}, we pick a proper subset $I
\subset \{ 1, \dots, s\}$ such that the union of the facets in $D_I$
is homeomorphic to a disk. In Section \ref{s:construct}, while
describing the algorithmic construction of the matrix of \ref{thm1},
we also show how to pick such $D_I$ as a partial shelling of the
facets of $Q$.  We will consider the sheaf $\F = \O(2D_Q - D_I)$. As
before this is a divisor on the normal toric variety $X_{\Sigma_Q}$
pushed forward onto $X_A$.  The main fact we will need is the
following cohomology vanishing theorem. For simplicity, and when there
is no confusion, we will often write $H^i(\O(D))$ instead of $H^i(X_A,
\O(D))$.

\begin{thm}
Let $X= X_Q$ be a projective toric variety of dimension $n$ arising
from a polytope $Q$ with corresponding ample divisor $D_Q$. Let $D_I$
be a proper subset of the facets such that the unions of the facets in
$D_I$ is a topological manifold with no reduced homology. Then:
\label{t:cohomvanish}
\begin{align*}
H^0 (\O(kD_Q - D_I) &\cong S_{kQ - D_I} & \\
H^i (\O(kD_Q - D_I) &\cong 0 & \quad  i = 1,\dots, n-1 \\
H^n (\O(kD_Q - D_I) &\cong S^{\ast}_{-kQ - \overline D_I}&
\end{align*}
\noindent for all $k \in \ZZ$.
\end{thm}

In the case $Q$ is a 3-polytope the only 2-manifold with no reduced
homology is the disk.  The proof is postponed until Section
\ref{ss:cohomology}.  But note that plugging this into the description
of the Tate resolution using $\F(k) = \O((k+2)D - D_I)$ gives us:

\begin{cor} \label{cor:tatedim3}
The Tate resolution of $\F$ has terms:
\begin{align*}
T^e(\F)     &\cong S^{\ast}_{(1-e) Q - {\overline {D_I}}}
 \otimes E(3-e) \quad {\rm for } \  e < -1 \\
T^{-1}(\F)  &\cong S^{\ast}_{2 Q - {\overline {D_I}}}
 \otimes E(4) \oplus S_{Q - D_I} \otimes E(1) \\
T^0(\F) &\cong S^{\ast}_{Q - {\overline {D_I}}}
 \otimes E(3) \oplus S_{2 Q - D_I} \otimes E \\
T^e(\F)  &\cong S_{(e + 2) Q - D_I} \otimes E(-e)
 \quad {\rm for } \  e > 0.
\end{align*}
\end{cor}

Finally, to get the Chow form we need to apply the functor $U_4$
which sends $E(p)$ to $\wedge^p \TT$.  But, $\TT$ is a vector
bundle of rank $4$, so by the above proposition only $T^{-1}(\F)$
and $T^0(\F)$ survive the application of $U_4$.  Therefore, $\F$
is weakly Ulrich and the matrix of the resulting two term complex
is exactly the matrix of Theorem \ref{thm1} which we restate here
in a slightly different language.

\begin{cor}
\label{cor:twotermdim3}

The resultant of $f_1, \dots, f_4$ is the determinant of the two
term complex below:

\begin{diagram}
&&S^\ast_{2 Q - \overline {D_I}} \otimes \bigwedge^4 \TT  &
                       \rTo^{\tilde{L}} & S^\ast_{Q - \overline {D_I}}
                       \otimes \bigwedge^3 \TT &      & \\
0& \rTo & \oplus  & \rdTo^B  & \oplus  & \rTo & 0 \\
&&S_{Q - D_I}  \otimes \bigwedge^1 \TT &
                        \rTo^{L} & S_{2 Q - D_I}
                        \otimes \bigwedge^0 \TT
\end{diagram}

\end{cor}

Theorem \ref{t:cohomvanish} can be used to give exact determinantal
formulas for resultants in dimension 4 and above for some cases of polytopes.

\begin{thm} \label{t:dim4res} 
Let $Q \subset \RR^4$ be a polytope such that $A = Q \cap \ZZ^4$
affinely spans $\ZZ^4$. There is a determinantal formula for $Res_A$
if $Q$ has no interior points and there is some facet $D_i$ of $Q$
with no relative interior points.
\end{thm}

\begin{pf} Take $\F = \O(2D_Q - D_i)$. Going through
the Tate resolution machinery using our vanishing theorem, we get a
three term complex whose left most term is
$S_{Q-\overline{D_i}}^{\ast}$. The points here are exactly the
interior points of $Q$ together with the relative interior points of
$D_i$. So under the given hypothesis, this term is zero and we
have a two term complex.
\end{pf}

In the case of $X_Q = \PP^4$ we recover the formulas for resultants of
5 homogeneous polynomials of degree less than or equal to 3.  We can
make a similar statement in dimension 5 and higher but the hypotheses
get stricter.

\begin{thm} \label{t:hidimres}
 Let $Q \subset \RR^n$ and $A = Q \cap \ZZ^n$ affinely spans $\ZZ^n$
for $n \geq 5$. let $k_1 = \floor{\frac{n+1}{2}} - 2$ and $k_2 =
\ceil{\frac{n+1}{2}} -2$. There is a determinantal formula for $\res_A$
if there is a collection of facets $D_I$ of $Q$ forming a manifold
without homology such that $k_1Q$ and $k_2Q$ have no interior lattice
points, $D_I$ has no relative interior lattice points in $k_2Q$ and
$\overline{D_I}$ has no relative interior points in $k_1Q$.
\end{thm}

\begin{pf} Take $\F = \O(\floor{\frac{n+1}{2}}Q - D_I)$. The result
follows from the Theorem \ref{t:cohomvanish} and counting lattice points.
\end{pf}
For example when $n = 5$ and $Q$ is the coordinate simplex we recover
the determinantal formula for $6$ homogeneous polynomials of degree
$2$. For $n=6$ or greater we only get a resultant formula for $d=1$.
It would be interesting to classify all polytopes of arbitrary shape
satisfying these conditions. It may be that there is only be a finite
list for $n=6$ or greater.

We do not claim that these theorems generate all determinantal
resultant formulas. For example, by Proposition 2.6 of \cite{ES} if
$Q_1$ and $Q_2$ (of any dimension) have resultant formulas with
sheaves $\F_1$ and $\F_2$, then $\F_1 \otimes \F_2$ will give a
determinantal formula for $Q_1 \times Q_2$. In any case polytopes
satisfying Theorem \ref{t:hidimres} together with all products of
such polynomials is at least a start towards classifying exact
resultant formulas in higher dimension.

\subsection{Cohomology vanishing}

\label{ss:cohomology}

In this section we will prove Theorem \ref{t:cohomvanish}.  So we will
need to compute the cohomology of $\O(kD_Q - D_I)$ for all $k \in
\ZZ$. We already know the global sections $H^0(X_A, \O(\cdot))$. The
next proposition shows how to compute the top cohomology $H^n(X_A, \O(\cdot))$.

\begin{prop} \label{p:topcohom}
Let $Q \subset \RR^n$ be a lattice polytope of dimension $n$ with
facets $D_1, \dots, D_s$ and $A = Q \cap \ZZ^n$ affinely
generating $\ZZ^n$. Let $X_A$ be the corresponding toric variety,
and $D = \sum a_i D_i$ a Weil divisor on the normalization
$X_{\Sigma_Q}$ which pushes forward as before to a sheaf on $X_A$.
Then
$$ H^n(X_A, \O(D)) \cong H^0(X_A, \O(-D - \sum_{i=1}^s D_i))^\ast $$
\end{prop}

\begin{pf}
As per our earlier discussion all of the cohomology can be
computed on the associated normal toric variety $X =
X_{\Sigma_Q}$. This is Cohen-Macaulay with dualizing sheaf
$\omega_{X} = \O(-\sum_{i=1}^s D_i)$.  If $D$ were Cartier the
statement would follow immediately from Serre duality. In the
general Weil divisor case we have to be a little bit more careful.
So we compute:

\begin{align*}
H^n(X, \O(D))^{\ast} &\cong \Hom(\O(D), \omega_{X}) \\
 &\cong \Hom(\O(D), \mathcal{H}om(\O(\sum_{i=1}^s D_i), \O_{X})) \\
 &\cong \Hom(\O(D) \otimes \O(\sum_{i=1}^s D_i), \O_{X}) \\
 &\cong \Hom(\O_{X}, (\O(D) \otimes \O(\sum_{i=1}^s D_i))^{\ast}) \\
 &\cong H^0(X, (\O(D) \otimes \O(\sum_{i=1}^s D_i))^{\ast}).
\end{align*}

The first isomorphism is Serre duality. The second uses that Weil
divisors are reflexive sheaves and $\mathcal{H}om(\O(D), \O_{X})
\cong \O(D)^{\ast} \cong \O(-D)$. The third and fourth steps are
by the adjointness of $\mathcal{H}om$ and $\otimes$, and the last
step is the definition of global sections. Finally, by Corollary
2.1 in \cite{Har}, the dual of any coherent sheaf is reflexive.
So,
$$(\O(D) \otimes \O(\sum_{i=1}^s D_I))^{\ast} \cong (\O(D) \otimes
\O(\sum_{i=1}^s D_I))^{\ast\ast\ast}.$$
\noindent However $(\O(D) \otimes \O(E))^{\ast\ast}$ is always
isomorphic to $\O(D+E)$ even if $D$ and $E$ are not locally free.
Hence we get
$$(\O(D) \otimes \O(\sum_{i=}^s D_I))^{\ast} \cong \O(-D - \sum_{i=1}^s D_i),$$\noindent as desired.
\end{pf}

It remains to show that the ``middle cohomology'' always vanishes
under the given conditions.   The proof is broken up into three parts, showing
$H^i(\O(kD_Q - D_I)) = 0$ when $k > 0$, $k = 0$, and $k < 0$.  The
first two follow fairly easily from results of Musta\c{t}\u{a}
\cite{Mus}, \cite{EMS}. The case $k <0$ requires more work and will be quite
interesting in its own right.

\begin{prop}\label{p:vanish1} Let $Q$ be a polytope and  $X_A$ the toric variety as in Proposition \ref{p:topcohom}. Let $D_I$ be the sum of any collection of facets as
before. $H^i(\O(kD_Q - D_I)) = 0$ for all $i > 0$ and all $k
> 0 $.
\end{prop}

\begin{pf}
Since $kD_Q$ is ample, this is just \cite[Corollary
2.5(iii)]{Mus}. 
\end{pf}

In general, the cohomology of all divisors can be grouped into a single object $H^i_{\ast}(\O_X)$ which has a $\ZZ^s$ fine grading:
$$H^i_{\ast}(\O_X) = \oplus_p H^i_{\ast}(\O_X)_p.$$
\noindent where $p \in \ZZ^s$.

The cohomology of a particular divisor class $[D]$ can now 
be recovered as 
$$H^i(\O_X(D)) = \sum_p H^i_{\ast}(\O_X)_p$$
\noindent where the sum is over all $p$ such that  $[\sum p_iD_i] = [D]$.

The next lemma can be viewed as a reformulation of a result of
\cite{Mus} yielding a topological formula for computing these graded
pieces. It shows that in the case of a projective toric variety sheaf
cohomology can be computed in terms of the ordinary homology of pure
cell complexes.
 
\begin{lem}< \label{l:topform} Let $p \in \ZZ^s$. Let $J =
{\rm neg}(p) \subset \{1, \dots, s\}$ be the set of coordinates for
which $p$ is strictly negative.  Let $|D_J|$ be the topological space
consisting of the union of all the facets $D_j$ with $j \in J$ of the
polytope $Q$ of $X$.
$$H_{\ast}^i(\O_X)_p \cong \tilde{H}^{i-1}(|D_J|).$$
\noindent The latter is the ordinary reduced cohomology of  $|D_J|$.
\end{lem}

\begin{pf}

Let $Y_J$ be the union of all cones in the fan $\Sigma$ having all
edges in the complement of $J$. Theorem 2.7 in \cite{EMS} shows
that for $i \geq 1$:
$$H_{\ast}^i(\O_X)_d \cong \tilde{H}^{i-1}(\RR^n \setminus Y_J).$$
The latter is isomorphic to $\tilde{H}^{i-1}(S^{n-1} \setminus
S^{n-1} \cap Y_J)$ (excision) which is further isomorphic to
$\tilde{H}_{n-i-1}(S^{n-1} \cap Y_J)$ by topological Alexander
duality. This is a subcomplex of the boundary complex of a
polytope polar dual to $Q$. The combinatorial Alexander dual is
the set of faces of $Q$ whose dual is not in $Y_I$. But this is
precisely all of those faces of $Q$ contained in some facet $D_I$.
The underlying topological space is $|D_I|$. So, by combinatorial
Alexander duality:
$$\tilde{H}_{n-i-1}(S^{n-1} \cap Y_J) \cong \tilde{H}^{i-1}(|D_J|).$$
\noindent as desired

\end{pf}

We now tackle the case $k = 0$, the proof of this next proposition was
given to me in a personal communication with Mircea Musta\c{t}\u{a}.

\begin{prop} \label{p:vanish2}
If the union of the collection of facets in $D_I$ is non-empty and
homologically trivial, then $H^i(\mathcal{O}(-D_I)) = 0$ for {\em all}
$i$. More generally, \\ $H^i(\O(-D_I)) \cong \tilde{H}^{i-1}(|D_I|)$,
\end{prop}

\begin{pf} (Due to Musta\c{t}\u{a})

$H^0(\O(-D_I)) = 0$ as the corresponding polytope is empty. 
Let $p_I$ be such that $(p_I)_i = -1$ if $i \in I$ and $(p_I)_i = 0$,
otherwise. Clearly, ${\rm neg}(p_I) = I$ and $\sum (p_I)_i D_i =
-D_I$.  By Lemma \ref{l:topform}, $H^i_{\ast}(\O_X)_{p_I} =
\tilde{H}^{i-1}(|D_I|)$.

 We now show that if $q$ is such that $[\sum q_i D_i] = [-D_I] $, but
$q \neq p_I$, then $H^i_{\ast}(\O_X)_q = 0$ for all $i$. Indeed, by
linear equivalence $q = p_I + {\rm div}(u)$, for some $u \in \ZZ^d \neq
0$. Let $J = {\rm neg}(q)$.  It is clear that
$$J = \{i|\langle u, \eta_i \rangle < 0 \ {\rm or} \ \langle u, \eta_i \rangle = 0 \ {\rm and} \ i \in I \}$$
Now the above implies there is a hyperplane $H \subset \mathbb{R}^s$
which separates the edges of $\Sigma_Q$ indexed by $J$ and
$\overline{J}$. By \cite[Proposition 2.6]{EMS} this forces
$H^i_{\ast}(\O_X)_q = 0$ for $i \geq 1$. 

\end{pf}

To complete the proof of Theorem \ref{t:cohomvanish} we need to consider
the case $k < 0$. This will require a new vanishing theorem which has
intrinsic interest.  Therefore, we state it in somewhat more
generality than necessary.

\begin{thm} Let $X$ be a projective toric variety of dimension
$n$, and $D$ a nef and big line bundle on $X$. Let $D_I = \sum_{i \in
I} D_i$ be a sum of prime torus invariant divisors. If the union of
the facets $D_i$ with $i \in I$ of $Q$ is a topological manifold with
boundary then $H^i(X, \O(-D-D_I)) = 0$ for all $0 \leq i < n$.
\end{thm}

Proposition 3.3 in \cite{Mus} states that the fan of $X$ refines the
normal fan of $P_D$ and $\O(D)$ is the pull-back of an ample divisor,
thus we can reduce to the case that $D$ is ample.

Theorem \ref{t:cohomvanish} gives general vanishing conditions for all
$k \in \ZZ$ but the results in this section show that the vanishing
theorem can be refined using different hypotheses for different cases
on the integer $k$.  When $k > 0$, all higher cohomology vanishes for
any subset $D_I$. When $k = 0$ we need the toplogical space $|D_I|$ to
have no reduced homology in which case all cohomology
vanishes. Finally for $ k <0$, when $|D_I|$ is a manifold, all
cohomology vanishes except at the top.

\begin{pf} 
By the remark above assume that $D$ is ample.  As before we will need
to compute $H^i_{\ast}(\O_X)_p$ for $\sum p_i D_i$ linearly equivalent
to $-D-D_I$. Let $p_I$ be defined as in the proof of Proposition
\ref{p:vanish2}. Any $p$ as above is of the form $q - p_I $ where
$\sum q_iD_i$ is linearly equivalent to $-D$.  Write $D = \sum a_i
D_i$, in which case $q_i = \langle u, \eta_i \rangle - a_i$ for some
$u \in \ZZ^n$.

Therefore, 
\begin{align*}
{\rm neg}(q) &= \{i|\langle u, \eta_i \rangle < a_i\} & \text{and}\\
{\rm neg}(p) &= {\rm neg}(q) \cup \{i | \langle u, \eta_i \rangle = a_i \ {\rm and} \ i \in I \} & 
\end{align*}
Let $J' = {\rm neg}(q)$ and $J = {\rm neg}(p)$ with $|D_{J'}|$ and
$|D_J|$ the corresponding unions of facets. Since $D$ is an ample
divisor, $H^i(\O(-D)) = 0$ for $i < n$, derived for example by
Proposition \ref{p:vanish1} and Serre duality. We need to show that
under the given hypotheses $H^i(\O(-D-D_I)) = 0$. We know by Lemma
\ref{l:topform} $\tilde{H}^i(|D_{J'}|) = 0$ for $i < n-1$, and so 
will need to prove that $\tilde{H}^i(|D_J|) = 0$.

We have three cases for $u$:

\begin{description}

\item[Case 1] $\langle u, \eta_i \rangle < a_i$ for all
$i$. Equivalently, $-u \in {\rm int}(P_D)$. In this case $|D_J|$ is the
entire boundary of $P_D$ which is an $n-1$ sphere and only
has reduced homology at the top.

\item[Case 2] $\langle u, \eta_i \rangle \leq a_i$ for all
$i$ and $\langle u, \eta_i \rangle = a_i$ for some $i$.

This means that $-u$ is on the boundary of $P_D$.  Since $D$ is ample,
$P_D$ has the same normal fan as $Q$ and so parallel faces to $Q$. The
set of all facets $D_j$ for which $\langle u, \eta_j \rangle = a_j$
cuts out a face $f$ of $Q$. Moreover, since $-D$ is Cartier there is a
corresponding function $\psi_{-D}$ on the fan $\Sigma$, defined to be
$a_i$ on the rays $\eta_i$ and extended linearly in each cone. Since
the linear functional $\langle u, \cdot \rangle$ agrees with
$\psi_{-D}$ on a spanning set of the cone corresponding to $f$ it
agrees with $\psi_{-D}$ on all of this cone.  Therefore, $\langle u,
\eta_i \rangle = a_i$ for all facets $D_i$ containing $f$ and so
$|D_J'|$ is the union of all facets of $Q$ not containing $f$. If $f$
is not a face of a facet in $D_I$ then none of the $D_j$ above are
part of $D_I$, in which case ${\rm neg}(p) = {\rm neg}(q)$ and
therefore $H^i(|D_J|) = H^i(|D_{J'}|)$.

Next, assume that $f$ is a face of some facet in $D_I$. The facets
$D_I$ define a cell complex, also denoted $D_I$, realizing the
manifold $|D_I|$.  The star $st(f)$ is the union of all of the relatively open faces of $D_I$ that have $f$ as a face and the link $lk(f)$ is $\overline{st(f)}- st(f)$. 

The key observation is that $|D_J| = |D_{J'}| \cup \overline{st(f)}$
and $lk(f) = |D_{J'}| \cap \overline{st(f)}$.  So, we have a
Mayer-Vietoris sequence:
$$ \cdots  \to \tilde{H}^{a-1}(lk(f))  \to \tilde{H}^a(|D_J|) \to  \tilde{H}^a(\overline{st(f)}) \oplus \tilde{H}^a(|D_{J'}|) \to \cdots.$$
We know that $\overline{st(f)}$ is contractible (it is star shaped!)
and from above $\tilde{H}^a(|D_{J'}|) = 0$ for $a < n-1$. It remains
to show that $\tilde{H}^{a-1}(lk(f)) = 0$ for $a < n-1$. This is where
we use that $|D_I|$ is a manifold.

Start with the cell complex $D_I$ and perform a stellar subdivision at
the face $f$. This induces a subdivision of $D_I$, which we call
$D_I^f$, with a new vertex $v_f$ corresponding to the face
$f$. Furthermore the star and link $st(v_f)$ and $lk(v_f)$ in $D_I^f$
are the same as $st(f)$ and $lk(f)$ in $D_I$. So it now suffices to
show that $\tilde{H}^{a-1}(lk(v_f)) = 0$ for $a < n-1$.

Since $|D_I|$ is a manifold with boundary, the local cohomology of
$|D_I|$ at $v_f$, $H_{v_f}^a(|D_I|)$, vanishes for $a \neq n-1$ if
$v_f$ is an interior point of $|D_I|$, and for all $a$ if $v_f$ is on
the boundary.  This local cohomology can also be computed from the
triangulation as the relative cohomology $H^i(\overline{st(v_f)},
lk(v_f))$.  The long exact sequence in relative cohomology yields:
$$ \cdots \to \tilde{H}^{a-1}(st(v_f)) \to \tilde{H}^{a-1}(lk(v_f)) \to
H^a(st(v_f), lk(v_f)) \to \cdots $$
Since $st(v_f)$ is contractible, $H^a(st(v_f), lk(v_f)) \cong
\tilde{H}^{a-1}(lk(v_f)) = 0$ for $a < n-1$ as desired.

\item[Case 3] $\langle u, \eta_i \rangle > a_i$ for some $i$.

In this case $-u$ is outside the polytope $P_D$.  A point $p$ in $P_D$
is {\em visible} from $-u$ if the straight line from $p$ to $-u$ meets
$P_D$ first in $p$. It is easy to see that if a visible point $p$ is
in the relative interior of a face $f$ then the whole face is visible
and any subface of a visible face is visible. Therefore visibility is
a property of whole faces. A face $f$ of $P_D$ will be called
{\em degenerate} if $-u$ is in the affine span of $f$. In particular $P_D$
itself is a degenerate face. A face is invisible if it is not visible
or degenerate. Any facet containing an invisible face must be
invisible or degenerate and if every facet containing some face $f$ is
degenerate then $f$ itself is degenerate.  Clearly, a facet $f$ is
visible if and only if $-u$ is on the opposite side of $f$ as
$P_D$. Therefore, $D_{J'}$ is the set of invisible facets. $D_J$ is
the union of $D_{J}'$ with some degenerate facets. So it will suffice to
prove the following proposition taking $P = P_D$ and $v = -u$.
\end{description}

\begin{prop} Let $P \subset \RR^n$ be a polytope of any dimension.
If $v$ is any point in the affine span of $P$ but outside of $P$, then
the union of the invisible facets of $P$ together with
any collection of degenerate faces is homologically trivial.
\end{prop}

\begin{pf} We proceed by induction on the number of degenerate faces and
the dimension of $P$. If $f$ is a degenerate face of $P$, so that $v$
is in the affine span of $f$ we can talk about the visible, invisible,
and degenerate faces of $f$ regarded as a polytope in its own right.
It is immediate from the definitions that a face of $f$ is visible
(invisible, degenerate) if and only if it visible (invisible,
degenerate) as a face of $P$. 

To apply the induction we need to show that the intersection of a
degenerate face $f$ with the union of the invisible facets and some
degenerate faces of $P$ is precisely the union of the invisible facets
of $f$ and some degenerate subfaces.

We first consider the intersection of a degenerate face $f$ with the
union of the invisible facets of $P$. Any invisible facet of $f$ is an
invisible face of $P$ and hence contained in an invisible facet of
$P$. For the converse, let $H$ be the affine span of $f$. Suppose $f'$
is a face of $f$ contained in an invisible facet $F$ of $P$. Since $u$
is on the same side of $F$ as $P$, it is on the same side of the
intersection of $F$ and $H$ as $f$. In particular there must be some
facet of $f$ containing $f'$ invisible from $u$. Hence, the union of
the invisible facets of $P$ intersects $f$ precisely in the union of
its invisible facets.

Next let $f$ be the intersection of two degenerate faces. Let $H$ be
the intersection of the corresponding two affine spans. So $H$
contains both $v$ and $f$ and moreover $H \cap P = f$. Let $H'$ be the
affine span of $f$, a subspace of $H$.  Each facet of $P$ defines a
half space containing $P$.  The intersection of all of these half
spaces for the facets containing $f$ is the convex hull of $P$ and $H'$.
Intersecting with $H$ yields just $H'$. One can instead take all of the
opposite half spaces and it remains true that the intersection with
$H$ is $H'$. Now if none of these facets are invisible from $v$, then
$v$ lies in all of the opposite half spaces as above, which means that
$v$ lies in $H'$ and thus $f$ is degenerate. In conclusion, the
intersection of two degenerate faces must either be degenerate or
contained in an invisible facet.

We can now proceed with the induction on the number of degenerate
faces. Let $P_0$ be the union of all the invisible facets of $P$.
This has no reduced cohomology since it is the negative support of a
negative ample divisor as before and therefore has no cohomology below
$\tilde{H}^{n-1}$.  Since $v$ is outside of $P$ there is at least one
visible facet and so the set of invisible facets is not the whole $n-1$-sphere.Therefore, $\tilde{H}^{n-1}$ is also 0.  

Assume now that $P_i$, the union of $P_0$ with $i$ degenerate faces,
is cohomologically trivial.  Let $f$ be a new degenerate face. The
Mayer-Vietoris sequence gives us:
$$ \cdots  \to \tilde{H}^{a-1}(f \cap P_i)  \to \tilde{H}^a(f \cup P_i) \to  \tilde{H}^a(f) \oplus \tilde{H}^a(P_i) \to \cdots.$$
As $f$ itself is contractible and $P_i$ is homologically trivial by
induction, it suffices to show that $f \cap P_i$ is homologically
trivial. However, the above arguments show that $v$ is in the affine
span of $f$ and $f \cap P_i$ is a union of all of the invisible facets
of $f$ and some degenerate faces of $f$. Therefore, its cohomology
vanishes by induction on dimension.  The base case is when $P$ is one
dimensional, in which case for $v$ in the line containing $P$ but not
in $P$, there is exactly one invisible facet ( a single point) and no
degenerate facets.
\end{pf}
Note, that this proposition, and hence all of Case 3, holds for
arbitrary $D_I$ and does not use that $D_I$ is a manifold.  Theorem
\ref{t:cohomvanish} is an easy consequence of all of the above
results.
\end{pf}

\section{Constructing the resultant matrix}
\label{s:construct}

\subsection{Partial shellings}
\label{ss:shelling}

In this section we show how to choose the $D_i$ to form a topological
ball (disk in dimension 2). Of course one can always choose a single
facet for $D_i$, but as we shall see this does not usually yield the
smallest matrices.

\begin{defn}
An ordering of the facets $D_1, \dots, D_s$ of an $n$-dimensional
polytope $Q$, is called a {\em shelling} if for $i= 2, \dots, s$, $(D_1
\cup \cdots \cup D_{i-1}) \cap D_i$ is $n-2$ dimensional and is itself
the union of an initial sequence of facets (codimension 2 faces in
$Q$) of a shelling of $D_i$. A {\em partial shelling} is a proper
sequence of facets, say $D_1, \dots, D_t$ with $1 \leq t < s$,
satisfying the same property above.
\end{defn}

When $Q$ has dimension $2$, a partial shelling is the same as a
connected set of edges.  In our setting, where $Q$ has dimension 3,
being a partial shelling simply means that the intersection of each
$D_i$ with the union of the previous $D_j$ is a connected set of edges
of $D_i$.  

\begin{prop} 
Let $Q$ be a polytope of dimension 3. The space $|D_I|$ is
homeomorphic to a disk if and only if the facets in $D_I$ can be
arranged into a partial shelling of the boundary of $Q$.
\end{prop}

\begin{pf} It is a standard result that any partial shelling
of the boundary of a polytope is homemorphic to a ball. In the
case of a 3 dimensional polytope it is actually a consequence of the
Jordan curve theorem.  Conversely, every topological disk is shellable
in dimension two. This last statement fails in dimension three and higher.
\end{pf}

It is very easy to actually construct partial shellings for
polytopes. A simple algorithm is to pass to the polar polytope
$Q^{\circ}$ of $Q$.  Facets of $Q$ correspond to vertices of
$Q^{\circ}$. Next, pick a generic vector in $\mathbb{R}^n$. This
will induce a linear functional on $Q^{\circ}$ which by genericity
induces a linear order on the vertices. One can show that any
initial segment of this linear ordering corresponds to a partial
shelling of the facets of $Q$. Shellings arising this way
are called line shellings. 

\subsection{Filling in entries}

\label{ss:entries}
 To actually construct our resultant formula we need to
fill in the entries of the matrices $B$, $L$, and $\tilde{L}$. We
saw above how these arise from a map in a Tate resolution.
Therefore, we must compute appropriate terms and maps in the Tate
resolution following the algorithm in Section \ref{tatealg} adapted to
this situation.

\begin{alg}

\smallskip \quad
\par
\begin{enumerate}
\item Pick a partial shelling $D_I$. As we shall see in the next
section, in order to get a smaller matrix we should pick $D_I$ to
have as many boundary points as possible.

\item Compute the lattice points in $3Q - D_I$ and $4Q - D_I$
respectively.

\item Construct the linear map $\phi_2 \ : \ S_{3Q-D_I} \otimes E \to
S_{4Q - D_I} \otimes E$. In light of Theorem \ref{cor:tatedim3} this
is precisely the differential $T^1(\F) \to T^2(\F)$.

Recall that $y_i$ represents a basis element of $S_Q$ hence a point in
$A = Q \cap \mathbb{Z}^3$. So, for every basis element $m$ of
$S_{3Q-D_I}$, let the multiplicative notation $m \cdot y_i$ denote the
basis element of $S_{4Q - D_I}$ obtained by adding the two
points. This can of course be extended linearly to all of
$S_{3Q-D_I}$. Now the map $\phi_2$ is explicitly defined by $\phi_2(m
\otimes 1) = \sum_{i=1}^N (my_i \otimes e_i)$.

\item Compute two steps of a graded minimal free resolution, over $E$,
of the cokernel of $\phi_2$.
$$ \begin{CD} T^{-1} @>\phi_0>> T^0 @>\phi_1>>T^1 @>\phi_2>> T^2
\end{CD} $$
Since this minimal free resolution is precisely the Tate resolution,
the map we are interested in is $\phi_0$. Let $M_0$ be the
corresponding matrix over $E$. The entries of this matrix will
be either linear or of degree $4$.

\item Apply the functor $U_4$ to $\phi_0$, and therefore $M_0$. This
is done by replacing each degree $4$ term of the form $e_{i_1}e_{i_2}e_{i_3}e_{i_4}$ by the ``bracket variable'' $[i_1i_2i_3i_4]$ which represents
the $4 \times 4$ determinant:
$$   \det \bmatrix  C_{1i_1} & C_{1i_2} & C_{1i_3} & C_{1i_4} \\
                C_{2i_1} & C_{2i_2} & C_{2i_3} & C_{2i_4} \\
                C_{3i_1} & C_{3i_2} & C_{3i_3} & C_{3i_4} \\
            C_{4i_1} & C_{4i_2} & C_{4i_3} & C_{4i_4} \endbmatrix$$
Here $C_{ij}$ is the coefficients of $f_i$ corresponding to the
monomial representing the point $y_j \in A$. These entries make up the
submatrix $B$ from Theorem \ref{thm1}. The remaining rows and columns
have linear entries, and correspond to $L$ and $\tilde{L}$. Replace
each such row (or column) by 4 rows (or columns).  The entry $e_i$ is
replaced by $C_{1i}$ in the first copy, $C_{2i}$ in the second copy,
and so on.  This procedure is illustrated in the examples below. It is
a consequence of Lemma $4.2$ in \cite{khe}.  This results in a matrix
$M$ which is precisely the matrix of Theorem \ref{thm1}.

\end{enumerate}
\end{alg}

Step 4 above requires computing part of a graded minimal resolution
over the exterior algebra. This can be done using Gr\"obner bases but
may be quite time consuming. On the other hand this computation needs
only be done once to compute the resultant of any system with a fixed
support.  One might hope to eliminate the expensive Gr\"obner basis
computations by finding explicit formulas for the non-trivial maps in
the resolution. This was done for the two-dimensional resultant in
\cite{khe} but remains open in the three-dimensional case.
\subsection{Examples}
\label{ss:examples}

\begin{figure}
\begin{center}
\epsfig{file=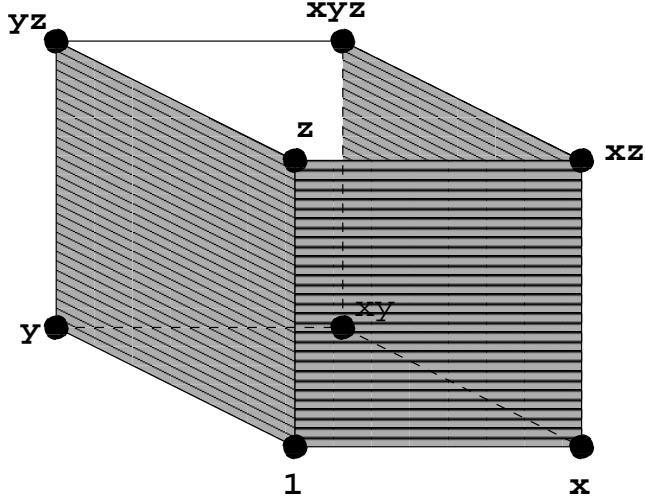} \caption{Newton Polytope of Example
\ref{ex:dim3_1}} \label{f:dim3_1}
\end{center}
\end{figure}

\begin{exmp}

\label{ex:dim3_1} Consider the multilinear system:

\begin{align*}
f_1 &= C_{11} + C_{12}x + C_{13}y + C_{14}z + C_{15}xy + C_{16}xz + C_{17}yz + C_{18}xyz \\
f_2 &= C_{21} + C_{22}x + C_{23}y + C_{24}z + C_{25}xy + C_{26}xz + C_{27}yz + C_{28}xyz \\
f_3 &= C_{31} + C_{32}x + C_{33}y + C_{34}z + C_{35}xy + C_{36}xz + C_{37}yz + C_{38}xyz \\
f_4 &= C_{41} + C_{42}x + C_{43}y + C_{44}z + C_{45}xy + C_{46}xz + C_{47}yz + C_{48}xyz \\
\end{align*}

The Newton polytope $Q$ of this system is the unit cube in Figure
\ref{f:dim3_1}. In order to apply the resultant algorithm we must
choose a partial shelling. So, for example, we can pick the left,
front, and, right faces as shown. Now $l(Q - D_I)$ and, by symmetry,
$l(Q - \overline{D_I})$ are empty while $l(2Q - D_I)$ consists of the 6
monomials $\{xy, xyz, xy^2, xyz^2, xy^2z, xy^2z^2\}$, while $l(2Q -
\overline{D_I})$ consists of the 6 monomials $\{z, xz, yz, x^z, xyz,
x^2yz\}$.  By Theorem \ref{thm1} the resultant is the determinant of a
$6\times 6$ pure B\'ezout matrix. To explicitly compute it, we
construct the linear map $S_{2Q-D_I} \otimes E(1) \to S_{3Q-D_I}
\otimes E$ and compute one step of a free resolution over $E$. The
matrix turns out to be the one shown in Table \ref{tbl:dim3_1}.

Note that the size of the matrix depends heavily on the choice of
the partial shelling. If, on the other hand, we were to choose
$D_I$ to consist of the left, front, and top facets, then $\#l(Q-
D_I) = \#l(Q-\overline{D_I}) = 1$, and $\#l(2Q-D_I) =
\#l(2Q-\overline{D_I}) = 8$. Hence, the matrix from Theorem
\ref{thm1} would be a $12\times 12$ matrix with an
$8\times 8$ block $B$, a $8\times 4$ block $L$, a $4\times 8$
block $\tilde{L}$, and a $4\times 4$ block of zeroes. 

If instead we tried the top and bottom facets, not homeomorphic to a
disk, we would still have $l(Q-D_I)$ and $l(Q-\overline{D_I})$
empty. However, this time $l(2Q-D_I)$ would consist of 9 points, while
$l(2Q-\overline{D_I})$ would have only 3 points. A closer look at the
vanishing theorems shows that we can still get a $9 \times 9$ square
resultant matrix as there is only other non-vanishing cohomology term
$H^1(\O(-D_I)) = \tilde{H}_0(|D_I|) = \CC$ tensored with $\bigwedge^2 \TT$,
a vector bundle of rank 6.  Indeed one can show in general that if
$|D_I|$ is a disjoint union of disks we still get an exact matrix
formula.

\end{exmp}

\begin{table}
\begin{center}
\caption{Resultant Matrix for Example \ref{ex:dim3_1}}
\label{tbl:dim3_1}
\end{center}
\begin{tiny}

\[ \bmatrix
 [1234]       & [1236]-[1245]& [1237]      & [1256]       & [1238]+[1257]& [1258]       \\
              &              &         &          &          &              \\
 [1346]-[1247]& [2346]-[1248]&[2347]-[1367]&-[2456]-[1268]& [2348]-[1368]&-[1568]-[2458]\\
          &-[1267]-[1456]&         &          &-[1567]-[2457]&              \\
 [1345]       & [2345]-[1356]&-[1357]      &-[2356]   &-[2357]-[1358]&-[2358]       \\
              &              &         &          &          &              \\
 [1467]       & [2467]+[1468]& [3467]      & [2468]   & [3468]-[4567]&-[4568]       \\
              &              &         &          &          &              \\
 [1457]+[1348]& [1368]+[2348]&[1378]-[3457]& [2368]-[2567]& [1578]+[2378]& [2578]-[3568]\\
          &+[2457]-[1567]&         &          &+[3458]-[3567]&              \\
-[1478]       &-[1678]-[2478]&-[3478]      &-[2678]       &
[4578]-[3678]& [5678]
 \endbmatrix
\]

\end{tiny}

\end{table}

\begin{figure}
\begin{center}
\epsfig{file=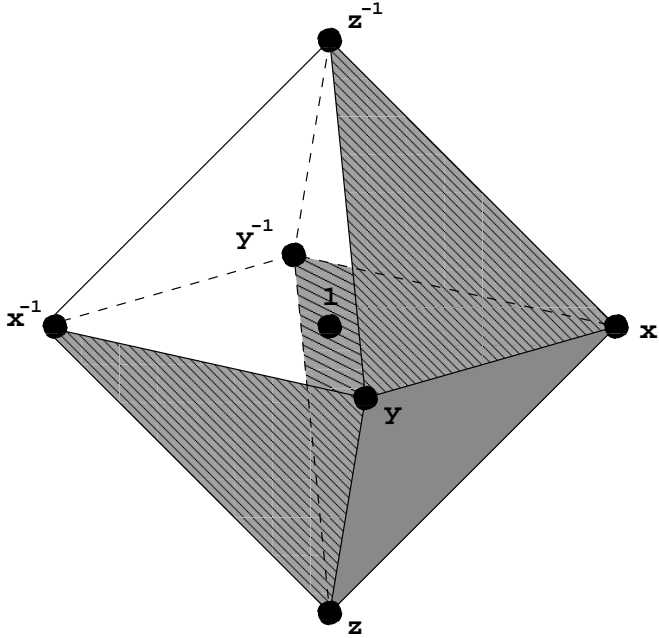} \caption{Newton Polytope of Example
\ref{ex:dim3_2}} \label{f:dim3_2}
\end{center}
\end{figure}

\begin{exmp}
\label{ex:dim3_2} Our next example is the following system:

\begin{align*}
f_1 &= C_{11} + C_{12}x + C_{13}y + C_{14}z + C_{15}x^{-1} + C_{16}y^{-1} + C_{17}z^{-1}\\
f_2 &= C_{21} + C_{22}x + C_{23}y + C_{24}z + C_{25}x^{-1} + C_{26}y^{-1} + C_{27}z^{-1}\\
f_3 &= C_{31} + C_{32}x + C_{33}y + C_{34}z + C_{35}x^{-1} + C_{36}y^{-1} + C_{37}z^{-1}\\
f_4 &= C_{41} + C_{42}x + C_{43}y + C_{44}z + C_{45}x^{-1} + C_{46}y^{-1} + C_{47}z^{-1}\\
\end{align*}

The Newton polytope $Q$ is the octahedron of Figure
\ref{f:dim3_2}. As our set of facets (partial shelling) we choose
the $x, y, z$ facet and the three other facets adjoined to it by
an edge. The chosen facets are shaded in the figure. Now we can
see that there are 10 points in $l(2Q - D_I)$ and also by symmetry
in $l(2Q - \overline{D_I})$. There is a single point in $l(Q - D_I)$
(respectively, $l(Q - \overline{D_I})$. By Theorem
\ref{thm1} the resultant is therefore the determinant of a
14x14 matrix shown in Table \ref{tbl:dim3_2}.  This matrix was
found following the algorithm of Section \ref{ss:entries} by
starting with the map $S_{3D_Q - D_I} \otimes E(1) \to S_{4D_Q -
D_I} \otimes E$ and computing a free resolution.

If we were to choose a non-partial shelling such as two facets meeting
at a single point then the corresponding resultant complex would have
nontrivial middle cohomology. Indeed, $H^2(-D-D_I) = H^2(-2D-D_I) =
\CC$.  $|D_I|$ is still homologically trivial but it is not a disk.
The complex arising from the Tate resolution has three terms.
\end{exmp}

\begin{table}

\begin{tiny}
\[\setcounter{MaxMatrixCols}{20}  \bmatrix
 0 &  0 & 0 & 0 & 0 & 0 & 0 & [2345] & [2346] & [2347] & C_{11} & C_{21} & C_{31} & C_{41} \\
 0 & 0  & 0 & 0 & 0 & 0 & 0 & 0  & 0 & 0 & C_{12} & C_{22} & C_{32} & C_{42} \\
 0 & 0  & 0 & 0 & 0 & 0 & 0 & 0  & 0 & 0 & C_{13} & C_{23} & C_{33} & C_{43}3 \\
 0 & 0  & 0 & 0 & 0 & 0 & 0 & 0  & 0 & 0 & C_{14} & C_{24} & C_{34} & C_{44} \\
 0 & -[2356]  & -[2357] & 0 & 0 & 0 & 0 & [1235]  & 0 & 0 & C_{15} & C_{25} & C_{35} & C_{45} \\
 0 & -[2456]  & 0 & [2467] & 0 & 0 & 0 & 0  & -[1246] & 0 & C_{16} & C_{26} & C_{36} & C_{46} \\
 0 & 0  & [3457] & [3467] & 0 & 0 & 0 & 0 & 0  & [1347] & C_{17} & C_{27} & C_{37} & C_{47} \\
 -[2567] & [1256]  & 0  & 0 & 0 & 0 & 0 &  -[2356] & -[2456] & 0  & 0 & 0 & 0 & 0 \\
 -[3567] & 0  & -[1357] & 0 & 0 & 0 & 0 & -[2357] & 0  & [3457] & 0 & 0 & 0 & 0 \\
 -[4567] & 0  & 0 & [1467] & 0 & 0 & 0 & 0 & [2467]  & [2467] & 0 & 0 & 0 & 0 \\
C_{11} & C_{12} & C_{13} & C_{14} & C_{15} & C_{16} & C_{17} & 0 & 0 & 0 & 0 & 0 & 0 & 0 \\
C_{21} & C_{22} & C_{23} & C_{24} & C_{25} & C_{26} & C_{27} & 0 & 0 & 0 & 0 & 0 & 0 & 0 \\
C_{31} & C_{32} & C_{33} & C_{34} & C_{35} & C_{36} & C_{37} & 0 & 0 & 0 & 0 & 0 & 0 & 0 \\
C_{41} & C_{42} & C_{43} & C_{44} & C_{45} & C_{46} & C_{47} & 0 &
0 & 0 & 0 & 0 & 0 & 0
 \endbmatrix
\]

\end{tiny}

\caption{Resultant Matrix for Example \ref{ex:dim3_2}}
\label{tbl:dim3_2}
\end{table}

\section{Ehrhart polynomials and sizes of resultant matrices}

\label{s:ehrhart}

The results of Section \ref{s:proof} show that the determinant of the
matrix of Theorem \ref{thm1} is the resultant. In particular, it must
be square and the degree of its determinant is equal to that of the
resultant. In this section we give an alternate combinatorial proof of
these facts. This will also allow us to analyze the size of the
resultant matrix in order to choose the smallest matrices.

Consider the Hilbert function of $X_A$, which turns out to be an
honest polynomial $p(x)$ where the value $p(k)$, for $ k \in
\mathbb{N}$, counts the number of lattice points in the polytope
$kQ$. This polynomial is associated to the polytope $Q$ and is
called the {\em Ehrhart polynomial} of $Q$.  There is a very pretty
duality theorem involving Ehrhart polynomials. See \cite{Ful} for
details.

\begin{prop}
Let $Q$ be a lattice polytope of dimension $n$ with Ehrhart
polynomial $p$. Then, $(-1)^{n}p(-k)$ is the number of {\em
interior} lattice points in $kQ$.
\end{prop}

Given a collection of facets $D_I$ we are interested in counting the
number of lattice points in $kQ - D_I$. A result of Stanley \cite{Sta}
extends Ehrhart polynomials and duality in this setting.

\begin{prop} \label{p:sta} \cite[Proposition 8.2]{Sta}
 
Let $Q$ be a lattice polytope of dimension $n$, and $D_I$ a collection
 of facets. Suppose $|D_I|$ is homeomorphic to a manifold. Then, there
 is a polynomial $p_I$ of degree $n$, such that $p_I(k)$ for $k>0$ is
 the number of points in $kQ - D_I$, $(-1)^{n}p_I(-k)$ for $k > 0$ is
 the number of lattice points in $kQ - \overline{D_I}$, and $p_I(0) =
 1 -\chi(|D_I|)$.  Here, $\chi(|D_I|)$ is the Euler characteristic of
 the manifold $|D_I|$. In particuular, if $|D_I|$ is a disk, then
 $p_I(0) = 0$.
\end{prop}

The difference $p(k) - p_I(k)$, the number of lattice points on 
the facets $D_I$ in $kQ$, is itself a polynomial of degree $(n-1)$.

Going back to resultants, we consider the two term
complex appearing in Corollary \ref{cor:twotermdim3}

\begin{diagram}
S^\ast_{2 Q - \overline {D_I}} \otimes \bigwedge^4 \mathcal{T}  &
                       \rTo^{\tilde{L}} & S^\ast_{Q - \overline {D_I}}
                       \otimes \bigwedge^3 \mathcal{T} &      & \\
\oplus  & \rdTo^B  & \oplus  & \rTo & 0 \\
S_{Q - D_I}  \otimes \bigwedge^1 \mathcal{T} &
                        \rTo^{L} & S_{2 Q - D_I}
                        \otimes \bigwedge^0 \mathcal{T}
\end{diagram}

Let $p_I(k)$ be the Ehrhart polynomial of $kQ - D_I$. This is a cubic
polynomial, thus the fourth difference is 0. In particular:
$$p_I(2) - 4p_I(1) + 6p_I(0) -4p_I(-1) + p_I(-2) = 0.$$
Since $|D_I|$ homeomorphic to a disk, $p_I(0) = 0$, and the equation
can be rewritten as $p_I(2) - 4p_I(-1) = -p_I(-2) + 4p_I(1)$.
Indentifying the dimension of the terms in the diagram above using
Proposition \ref{p:sta}, this says precisely that the matrix is square.

The total degree is computed by taking $4\#l(Q - D_I)$ entries from
$L$, $4\#l(Q-D_{\overline{I}})$ entries from $\tilde{L}$ and $\#l(2Q
- D_I) - 4\#l(Q-D_I)$ entries from $B$. The entries of $L$ and
$\tilde{L}$ are of degree 1, while those of $B$ are of degree $4$.
So the total degree is $4p_I(1) - 4p_I(-1) + 4(p_I(2) - 4p_I(1)) = 4(p_I(2)
- 3p_I(1) + 3p_I(0) - p_I(-1))$. This is $4$ times the third difference
of $p_I$ which is the same as $4$ times $3!$ times the leading
coefficient of $p_I$. This is the same as the leading coefficient of
the Ehrhart poynomial of $Q$ which is just the Euclidean volume.
Hence, the degree in question is $4$ times the normalized volume
which is also the degree of the resultant.

This leads to a technique to analyze the size of the resultant
marices.  The Ehrhart polynomial of $Q$ is of the form $p(x) = Ax^3
+ Bx^2 + Cx + 1$. The leading term $A$ is the degree of the toric
variety $X_A$ divided by $3!$, which is the Euclidian volume
of $Q$. Moreover $p(1)$ is the number of lattice points in $Q$,
and $p(-1)$ is the negative of the number of interior points. So
the number of boundary points in $Q$ is $p(1) + p(-1) = A + B + C
+ 1 -A + B -C +1 = 2B+2$. Let $B_Q = 2B + 2$ denote this number.
Next, for any partial shelling $D_I$, we write the corresponding
quadratic Ehrhart polynomial as $q_I(x) = ax^2 + bx + 1$. This time
$q(-1)$ is equal to the number of relative interior points, so the
number of boundary points is $q(1) - q(-1) = a + b + 1 -a + b -1 =
2b$. Let $B_I = 2b$ denote this number. Taking $p_I(x) = p(x) -
q_I(x)$ as above, then the total size of the resultant matrix is

\begin{equation*}
\begin{split}
p_I(2) - 4p_I(-1) &= p(2) - q_I(2) - 4p(-1) + 4q_I(-1) \\
              &= 8A +  4B + 2C + 1 - (4a + 2b + 1) - 4(-A + B - C + 1)\\
              & \quad  + 4(a - b + 1) \\
              &= 12A + 6C -6b. 
\end{split}
\end{equation*}

Let $i_Q = -p(-1) = A - B + C - 1$ be the total number of interior
points of $Q$. So $C = i_Q + 1 + B - A$. Hence, we can rewrite
the above as:

\begin{align*}
12A + 6C - 6b &= 12A + 6 (i_Q + B - A + 1)  - 6b \\ 
              &= 6A +  3(2B + 2 - 2b) + 6i_Q \\
              &= V + 3(B_Q - B_I) + 6i_Q. \\
\end{align*}

Here, $V$ denotes the normalized volume of $Q$ which is 6 times the
Euclidian volume $A$.  Therefore, in order to minimize the size of
the matrix we must maximize $B_I$ which is the number of relative
boundary points of the union of the facets $D_I$. This gives
an obvious lower bound of $V$ for the size of the resultant matrix.
A more sophisticated argument would give an upper bound of $3V$ 
when $i_Q$ is at least 1.

\section{Conclusion}

In this article we showed how the resultant of an unmixed system
in three variables with arbitrary support can be computed as the
determinant of a matrix. Combined with the authors earlier results
\cite{khe}, we have now generalized the formulas for the resultant
of homogeneous systems in dimensions 2 and 3. However, it is still
unknown how to make the dimension 3 formula completely explicit
instead of in terms of a free resolution as presented here.

For dimension 4 and higher, no general exact formula is known.  We do
give some special cases, although still without an explicit closed
form formula, and it would be nice to finish this classification.  A
second approach is to allow complexes with more than two terms,
yielding resultant formulas with extraneous factors. In the case of
projective spaces \cite{DiD} and products of projective spaces
\cite{DiE} the extraneous factors have been identified. It is still
open how to do this for general toric varieties.

\bibliographystyle{plain}
\bibliography{dim3res}

\end{document}